\documentclass[amsmath,amssymb,twoside,10pt]{article}
 \topmargin 0.3cm
\usepackage{amsfonts, amssymb}
\usepackage{graphicx,subfigure}

\usepackage{amsmath, amsfonts, amssymb, amsthm}
\newtheorem{theorem}{Theorem}

\newtheorem{lemma}{Lemma}[section]
\begin{document}
%
\setcounter{page}{27}  
%
\begin{center}
{\Huge INTERNATIONAL\\
\vspace{0.2cm}
 PUBLICATIONS USA}\\
\end{center}
\begin{center}
{\Large PanAmerican Mathematical Journal}\\
\vspace{0.1cm}
{\Large Volume 23(2013), Number 3, 27--34}\\
\vspace{0.4cm} {\bf{\large  Cycles in Random Bipartite Graphs}}
\\
\vspace{0.3cm}
Yilun Shang\\
University of Texas at San Antonio\\
Institute for Cyber Security\\
San Antonio, Texas 78249, USA\\
\tt{shylmath@hotmail.com}\\
\vspace{0.4cm} Communicated by Allan Peterson \\
{\scriptsize(Received January 2013; Accepted April 2013)}\\
\end {center}
\begin{abstract}\noindent
In this paper we study cycles in random bipartite graph $G(n,n,p)$.
We prove that if $p\gg n^{-2/3}$, then $G(n,n,p)$ a.a.s. satisfies
the following. Every subgraph $G'\subset G(n,n,p)$ with more than
$(1+o(1))n^2p/2$ edges contains a cycle of length $t$ for all even
$t\in[4,(1+o(1))n/30]$. Our theorem complements a previous result on
bipancyclicity, and is closely related to a recent work of Lee and
Samotij.
\\

\noindent {\bf AMS (MOS) Subject Classification:}   05C80
(05C38, 05D40) \\

\noindent {\bf Key words:} Random graph, bipartite graph, cycle.
\end{abstract}

\section{Introduction}
\newcommand{\Reals}{\mathbb R}
\newcommand{\Ts}{\mathbb T}
 \newcommand{\xproof}[1]{\par\noindent{\em Proof.  }\ #1 \, \nobreak \hfill $\blacksquare$ \vspace{.3pc}}
  \newcommand{\yproof}[1]{\par\noindent{\em Proof of Theorem 4.  }\ #1 \, \nobreak \hfill $\blacksquare$ \vspace{.3pc}}
 \newtheorem{defn}[theorem]{Definition}  
 \newtheorem{example}[theorem]{Example}
 \newtheorem{exer}{Exercise}
 \newtheorem{nonl}[theorem]{Lemma}  
 \newtheorem{non}[theorem]{Theorem} 
 \newtheorem{nonc}[theorem]{Corollary} 
 \newtheorem{nonp}[theorem]{Proposition}  

Given a complete bipartite graph $K_{n,n}$ and a real $p\in[0,1]$,
let random bipartite graph model $G(n,n,p)$ be the probability space
of subgraphs of $K_{n,n}$ obtained by taking each edge independently
with probability $p$ (see e.g. \cite{20}). For a given graph
property $\mathcal{P}$, we say that $G(n,n,p)$ possesses property
$\mathcal{P}$ asymptotically almost surely, or a.a.s. for brevity,
if the probability that $G(n,n,p)$ possesses $\mathcal{P}$ tends to
1 as $n$ goes to infinity. In the previous work \cite{1}, we
provided an edge condition for cycles in Hamiltonian subgraphs of
$G(n,n,p)$:

\begin{theorem}
If $p\gg n^{-2/3}$, then $G(n,n,p)$ a.a.s. satisfies the
following. Every Hamiltonian subgraph $G'\subset G(n,n,p)$ with more
than $(1+o(1))n^2p/2$ edges is bipancyclic (i.e., contains cycles of
every possible even length).
\end{theorem}

Bipancyclicity of bipartite graphs is first studied by Schmeichel
and Mitchem \cite{2}. Theorem 1 can be viewed as an extension in
random graph setting of a classical theorem of them \cite{3}, which
says that every Hamiltonian bipartite balanced graph with $2n$
vertices and more than $n^2/2$ edges is bipancyclic. We note that
Theorem 1 is best possible in two ways \cite{1}. First, the range of
$p$ is asymptotically tight. Second, the proportion $1/2$ of edges
cannot be reduced .

In the present work, we will focus on the situation where the
subgraph $G'$ is not necessarily Hamiltonian. We establish the
following result.

\begin{theorem}
If $p\gg n^{-2/3}$, then $G(n,n,p)$ a.a.s. satisfies the following.
Every subgraph $G'\subset G(n,n,p)$ with more than $(1+o(1))n^2p/2$
edges contains a cycle of length $t$ for all even
$t\in[4,(1+o(1))n/30]$.
\end{theorem}

Seeking cycles of various lengths in random graph models is an
interesting topic in probabilistic combinatorics and has since long
attracted much research attention, see e.g.
\cite{10,6,4,7,8,5,9,11}. Recently, Lee and Samotij \cite{12} showed
that if $p\gg n^{-1/2}$, then binomial random graph $G(n,p)$ a.a.s.
satisfies the following: Every Hamiltonian subgraph $G'\subset
G(n,p)$ with more than $(1/2+o(1))n^2p/2$ edges is pancyclic (i.e.,
contains cycles of every possible length). In the arXiv version of
their paper, they further proved the following result concerning
cycles of short length.

\begin{theorem}
For $\varepsilon\in(0,1)$, there exists a constant $C$ such that if
$p\ge Cn^{-1/2}$, then $G(n,p)$ a.a.s. satisfies the following.
Every subgraph $G'\subset G(n,p)$ with more than
$(1/2+\varepsilon)n^2p/2$ edges contains a cycle of length $t$ for
all $3\le t\le \varepsilon n/2560$.
\end{theorem}

To prove our main result, we will roughly follow the line of the
proof of Theorem 3. The bipartite structure of $G(n,n,p)$ entails
some significant modifications which remarkably allow us to show the
existence of longer cycles (i.e., cycles of linear length) with even
smaller edge probability (See Theorem 2).

The rest of the paper is organized as follows. We will present some
useful lemmas in Section 2 and prove the main result in Section 3.

\section{Some lemmas}

We begin with some notations. Let $G=(V_0,V_1,E)$ denote a bipartite
graph with two classes of bipartition $V_0$, $V_1$ and edge set $E$.
For a vertex $v$, we denote its $k$-th order neighborhood by
$N^{(k)}(v)$, i.e., the set of vertices at distance $k$ from $v$.
Let $N(v)=N^{(1)}(v)$ and $\operatorname{deg}(v)=|N(v)|$ be the
degree of $v$. The degree of set $X$ is defined as
$\operatorname{deg}(X)=\sum_{v\in X}\operatorname{deg}(v)$. The
first order neighborhood of set $X$ is defined as $N(X)=\{u\in
V_0\cup V_1:u\in N(v)\ \mathrm{for}\ \mathrm{some}\ v\in X\}$. The
maximum degree of graph $G$ is defined as $\Delta(G)$. For a set
$X$, let $E(X)$ be the set of edges in the induced subgraph $G[X]$,
and let $e(X)=|E(X)|$. Analogously, the set of ordered pairs
$(x,y)\in E$ with $x\in X$ and $y\in Y$ is denoted by $E(X,Y)$. Let
$e(X,Y)=|E(X,Y)|$. When there are several graphs under
consideration, we may use subscripts such as $N_G(v)$ to indicate
the graph we are currently working with. Floor and ceiling signs are
often omitted whenever they are not crucial.

The following concentration inequality (see e.g. \cite[Corollary
2.3]{13}) will be often used in the proof of main result.

\begin{lemma} (Chernoff's bound) Let $0<\varepsilon\le3/2$.
If $X$ is a binomial random variable with parameter $n$ and $p$,
then
$$
P(|X-\mathbb{E}(X)|\ge\varepsilon\mathbb{E}(X))\le2e^{-\varepsilon^2\mathbb{E}(X)/3},
$$
where $\mathbb{E}$ represents the expectation operator.
\end{lemma}

Lemma 2.2 can be found in \cite[Proposition 1.2.2]{14} and Lemma 2.3
is referred to as P\'osa's rotation-extension lemma (\cite{s} and
\cite[Chapter 10, Problem 20]{15}).

\begin{lemma} Let $G$ be a graph on $n$ vertices with at least $dn$ edges.
Then $G$ contains a subgraph $G'\subset G$ with minimum degree at
least $d$.
\end{lemma}

\begin{lemma}
Let $G=(V,E)$ be a graph such that $|N(X)\backslash X|\ge 2|X|-1$
for all $X\subset V$ with $|X|\ge t$. Then for any vertex $v\in V$,
there exists a path of length $3t-2$ in $G$ which has $v$ as an end
point.
\end{lemma}

For a monotone increasing property $\mathcal{P}$, the global
resilience \cite{16} of a graph $G$ with respect to $\mathcal{P}$ is
defined as the minimum number $r$ such that by deleting $r$ edges
from $G$, one can obtain a graph not possessing $\mathcal{P}$. The
following result was stated in \cite{1}, and can be proved similarly
as \cite[Proposition 3.1]{8} and \cite[Proposition 2.7]{12}.

\begin{lemma}
Assume that $0<p'\le p\le1$ and $n^2p'\rightarrow\infty$ as
$n\rightarrow\infty$. If $G(n,n,p')$ a.a.s. has global resilience at
least $(1/2-\varepsilon/4)n^2p'$ with respect to a monotone
increasing graph property, then $G(n,n,p)$ a.a.s. has global
resilience at least $(1/2-\varepsilon/2)n^2p$ with respect to the
same property.
\end{lemma}

In the next lemma we establish the expansion property for subgraphs
of $G(n,n,p)$ with large minimum degree. A variant of this result
for binomial random graph $G(n,p)$ appeared in \cite[Lemma 3.4]{8}.

\begin{lemma}
If $p=Cn^{-2/3}$ for some $C>0$ and $\varepsilon'\in(0,1)$, then
a.a.s. every subgraph $G'\subset G(n,n,p)$ with minimum degree at
least $\varepsilon'np$ satisfies the following expansion property.
For all $X\subset V_0\cup V_1$ with $|X|\le\varepsilon'n/15$, we
have $|N_{G'}(X)\backslash X|\ge2|X|$.
\end{lemma}

\xproof{Fix a subgraph $G'\subset G(n,n,p)$ with minimum degree at
least $\varepsilon'np$. Assume to the contrary that there exists
$X\subset V_0\cup V_1$ such that $|X|\le\varepsilon'n/15$ but
$|N_{G'}(X)\backslash X|<2|X|$. Let $Y=X\cup N_{G'}(X)$. We have
$|Y|\le 3|X|\le\varepsilon'n/5$. Since $G'$ has minimum degree at
least $\varepsilon'np$,
$$
e_{G'}(Y)=\frac12e_{G'}(Y,Y)\ge\frac12e_{G'}(X,X)\ge\frac12|X|\varepsilon'np\ge\frac16|Y|\varepsilon'np.
$$
Denote $|Y|=a$. Then we have $\varepsilon'np\le a\le
\varepsilon'n/5$.

The probability that there exists a set of order $a$ which spans at
least $a\varepsilon'np/6$ edges is
\begin{eqnarray}
&&\sum_{b=1}^{a-1}{n\choose b}{n\choose a-b}{(a-b)b\choose
a\varepsilon'np/6}p^{a\varepsilon'np/6}\nonumber\\
&\le&\sum_{b=1}^{a-1}\left(\frac{en}{b}\right)^b\left(\frac{en}{a-b}\right)^{a-b}\left(\frac{6e(a-b)b}{a\varepsilon'np}\right)^{a\varepsilon'np/6}
p^{a\varepsilon'np/6}\nonumber\\
&\le&\sum_{b=1}^{a-1}\frac{(en)^a}{b^b(a-b)^{a-b}}\left(\frac{3ea}{2\varepsilon'n}\right)^{a\varepsilon'np/6}.\label{1}
\end{eqnarray}
An application of Young's inequality yields
$$
\left(\frac1b\right)^b\left(\frac{1}{a-b}\right)^{a-b}\le\frac
ba\left(\frac1b\right)^a+\left(\frac{a-b}{a}\right)\left(\frac{1}{a-b}\right)^a.
$$
Therefore, the right-hand side of (\ref{1}) is at most
\begin{eqnarray*}
&&(en)^a\left(\frac{3ea}{2\varepsilon'n}\right)^{a\varepsilon'np/6}\cdot\sum_{b=1}^{a-1}\left(\frac{b^{1-a}}{a}
+\left(\frac{1}{a-b}\right)^a\right)\nonumber\\
&\le&(1+a)(en)^a\left(\frac{3ea}{2\varepsilon'n}\right)^{a\varepsilon'np/6}\ll\left(\frac{e}{3}\right)^{a\varepsilon'np/6}\ll
n^{-2}.
\end{eqnarray*}
Summing over all $\varepsilon'np\le a\le \varepsilon'n/5$, we see
that the probability that there is a set violating the assertion of
the lemma is $o(1)$. The proof is complete.}

\section{Existence of cycles}

In this section we will prove the following main theorem.

\begin{non}
For any $\varepsilon\in(0,1)$, there exists a constant $C$ such that
if $p\ge Cn^{-2/3}$, then $G(n,n,p)$ a.a.s. satisfies the following.
Every subgraph $G'\subset G(n,n,p)$ with more than
$(1+\varepsilon)n^2p/2$ edges contains a cycle of length $t$ for all
even $t\in[4,(1+\varepsilon/2)n/30]$.
\end{non}

The following proposition is a key step towards the proof of Theorem
4. Different from \cite{12}, the second order neighborhood of any
vertex in $G(n,n,p)$ contains no edges. We will resort to a coupon
collector argument to show the existence of a vertex with many edges
between its first order and second order neighborhoods.

\begin{nonp}
For any $\varepsilon\in(0,2/5]$, there exists a constant $C_0$
such that the following holds. If $p=Cn^{-2/3}$ for some constant
$C\ge C_0$, then $G(n,n,p)$ a.a.s. satisfies the following. Every
subgraph $G'\subset G(n,n,p)$ with more than $(1+\varepsilon)n^2p/2$
edges contains a vertex $v_0$ such that
$e(N_{G'}(v_0),B)\ge\varepsilon(1+\varepsilon/2)n^3p^3/8$ for some
$B\subset N_{G'}^{(2)}(v_0)$ satisfying $|B\cup
N_{G'}(v_0)|\le\varepsilon n^2p^2/2$.
\end{nonp}

\xproof{ Let $G$ be a graph drawn from $G(n,n,p)$ and $G'$ be a
subgraph with $e(G')>(1+\varepsilon)n^2p/2$. By Lemma 2.1, we have
a.a.s. $\Delta(G')\le\Delta(G)\le(1+\varepsilon)np$. In what
follows, we will show the proposition conditioned on the above
event.

For $i=0,1$, denote by $B_i$ the collection of the vertices in $V_i$
which have degree at least $(1+\varepsilon/2)np/2$ in $G'$.
Therefore,
\begin{eqnarray*}
(1+\varepsilon)n^2p\le 2e(G')=\operatorname{deg}_{G'}(V_0\cup
V_1)&=&\deg_{G'}(B_0\cup B_1)+\deg_{G'}((V_0\cup
V_1)\backslash(B_0\cup B_1))\\
&\le&\deg_{G'}(B_0\cup B_1)+(1+\varepsilon/2)n^2p,
\end{eqnarray*}
and hence
$$
\varepsilon n^2p/2\le\deg_{G'}(B_0\cup B_1)=\sum_{v\in V_0\cup
V_1}|N_{G'}(v)\cap (B_0\cup B_1)|.
$$
A simple averaging argument implies that there exists a vertex, say
$v_0\in V_0$, such that $|N_{G'}(v_0)\cap B_1|\ge\varepsilon np/4$
(c.f. Fig. 1).

\begin{figure}[h]
\begin{center}
\scalebox{0.6}{\includegraphics[119pt,224pt][580pt,570pt]{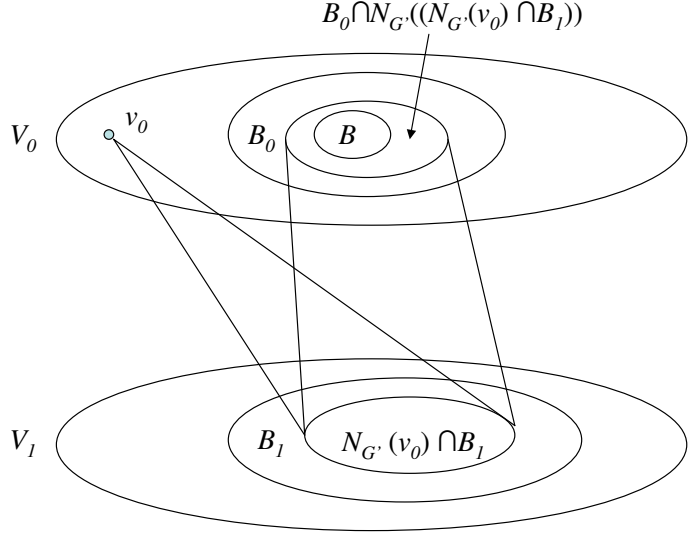}}\caption{A
depiction of bipartite graph $G$.}
\end{center}
\end{figure}

Consider a coupon collector's problem where a customer tries, in
several attempts, to collect a complete set of $|B_0|$ different
coupons. At each attempt, the collector gets a coupon randomly
chosen from $|B_0|$ kinds. We observe that the number $|B_0\cap
N_{G'}(N_{G'}(v_0)\cap B_1)|$ dominates the number $A$ of different
coupons the collector obtains in $T=\alpha|B_0|$ attempts, where
$$
\alpha=\frac{1}{|B_0|}\left(\frac{\varepsilon}{4}np\left(1+\frac{\varepsilon}{2}\right)\frac{np}{2}\right)=
\frac{1}{8|B_0|}\left(1+\frac{\varepsilon}{2}\right)\varepsilon
n^2p^2.
$$
Note that $|B_0|$ satisfies
$|B_0|(1+\varepsilon)np+(n-|B_0|)(1+\varepsilon/2)np/2\ge(1+\varepsilon)n^2p/2$,
which implies that $|B_0|\ge(\varepsilon/(2+3\varepsilon))n$.
Therefore,
$$
\alpha\le\frac{\left(1+\frac{\varepsilon}{2}\right)(2+3\varepsilon)}{8\varepsilon}np^2<1,
$$
for large enough $n$. Using Lemma 2.1, we can show that a.a.s. (see
e.g. \cite{17})
\begin{eqnarray*}
A&\ge&\left(1-\frac{\varepsilon}{2}\right)\left(1-\frac18\left(1+\frac{\varepsilon}{2}\right)\frac{\varepsilon
n^2p^2}{|B_0|}\right)\frac18\left(1+\frac{\varepsilon}{2}\right)\varepsilon
n^2p^2\\
&\ge&\frac18\left(1-\frac{\varepsilon}{2}\right)\left(1+\frac{\varepsilon}{2}\right)\varepsilon
n^2p^2\\
&>&\frac{\varepsilon}{4}n^2p^2,
\end{eqnarray*}
where the last inequality holds since $\varepsilon\le2/5$. According
to our above comment, we have
$$
|B_0\cap N_{G'}(N_{G'}(v_0)\cap B_1)|>\frac{\varepsilon}{4}n^2p^2.
$$

Now we choose a set $B\subseteq B_0\cap N_{G'}(N_{G'}(v_0)\cap
B_1)\subset N_{G'}^{(2)}(v_0)$ such that $|B|=\varepsilon n^2p^2/4$.
We obtain
$$
e(N_{G'}(v_0),B)\ge\left(1+\frac{\varepsilon}{2}\right)\frac{np}{2}|B|=\frac{\varepsilon}{8}\left(1+\frac{\varepsilon}{2}\right)n^3p^3,
$$
and
$$
|B\cup
N_{G'}(v_0)|\le(1+\varepsilon)np+\frac{\varepsilon}{4}n^2p^2\le\frac{\varepsilon}{2}n^2p^2
$$
as desired.}

Putting these together, we are now ready to show the main result.

\yproof{ Lemma 2.4 implies that it suffices to show the theorem with
$p=Cn^{-2/3}$ for some $C$ to be determined later. Without loss of
generality, we may also assume that $\varepsilon\le2/5$.  Let $G$ be
a graph drawn from $G(n,n,p)$ and $G'$ be a subgraph with
$e(G')>(1+\varepsilon)n^2p/2$.

It follows from Proposition 5, there exists a vertex, say $v_0\in
V_0$, such that
$e(N_{G'}(v_0),B)\ge\varepsilon(1+\varepsilon/2)n^3p^3/8$ for some
$B\subset N_{G'}^{(2)}(v_0)$ satisfying $|B\cup
N_{G'}(v_0)|\le\varepsilon n^2p^2/2$. Consider the subgraph
$G'[B\cup N_{G'}(v_0)]$. By Lemma 2.2, there exists a subset
$D\subset B\cup N_{G'}(v_0)$ such that $G'[D]$ has minimum degree at
least $(1+\varepsilon/2)np/4$. Take an arbitrary vertex $v_1\in
D\cap V_1$, and take a vertex $v_2\in D\cap V_0$ so that $v_0v_1v_2$
is a path of length 2 in $G'$. By Lemma 2.5, for all $X\subset
V_0\cup V_1$ of order $|X|\le(1+\varepsilon/2)n/60$, we have
$|N_{G'[D]}(X)\backslash X|\ge2|X|$. By Lemma 2.3, we can find a
path of length $(1+\varepsilon/2)n/30$ in $G'[D]$ which has $v_2$ as
an end point. We call this path $v_2x_1x_2x_3\cdots
x_{(1+\varepsilon/2)n/30}$, where $x_i\in N_{G'}(v_0)$ if $i$ is
odd; $x_i\in B\subset N_{G'}^{(2)}(v_0)$ if $i$ is even.

\begin{figure}[h]
\begin{center}
\scalebox{0.7}{\includegraphics[107pt,405pt][532pt,529pt]{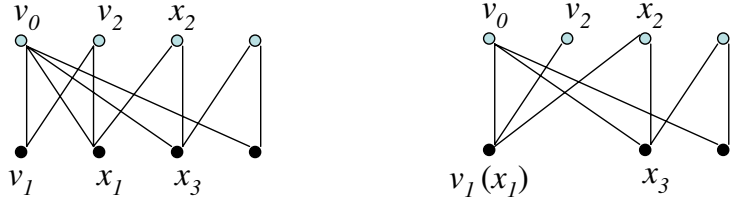}}\caption{A
depiction of case (i) and case (ii).}
\end{center}
\end{figure}

Now we consider two cases (c.f. Fig. 2):
\begin{itemize}
\item[(i)] $x_1\not=v_1$. Then $v_0v_1v_2x_1x_2x_3\cdots x_sv_0$ forms a
cycle of length $s+3$ in $G'$ when $s$ is odd;

\item[(ii)] $x_1=v_1$. Then $v_0x_1x_2x_3\cdots x_sv_0$ forms a cycle of
length $s+1$ in $G'$ when $s\ge3$ and $s$ is odd.
\end{itemize}

This completes the proof of Theorem 4.}


\end{document}